\documentclass[review]{elsarticle}
\usepackage{lineno,hyperref}
\usepackage{amssymb,amsmath,latexsym}
\usepackage{mathrsfs}
\usepackage{tipa}
\usepackage{amsfonts}
\usepackage{graphicx}
\usepackage{subfigure}
\usepackage[title]{appendix}
\usepackage{color}
\usepackage{multirow}
\usepackage{amsbsy}
\usepackage{longtable}
\usepackage{floatrow}
\usepackage{subfigure}
\usepackage{indentfirst}

\usepackage{ifpdf}

\newtheorem{theorem}{Theorem}
\newtheorem{definition}{Definition}
\newtheorem{proposition}{Proposition}
\newtheorem{lemma}{Lemma}
\newtheorem{cor}{Corollary}
\usepackage[font=small]{caption}
\modulolinenumbers[5]
\usepackage{xcolor}
\usepackage{hyperref} 
\hypersetup{               
	colorlinks=true,                
}

\AtBeginDocument{\hypersetup{pdfborder={0 0 1}}}

\usepackage[margin=2.5cm]{geometry}
\journal{}

\begin{document}
	
	\begin{frontmatter}
		
		\title{A note on sharp oracle bounds for Slope and Lasso}
		
		
		\author[mymainaddress1]{Zhiyong Zhou\corref{mycorrespondingauthor}}
		\cortext[mycorrespondingauthor]{Corresponding author}
		\ead{zhiyongzhou@zucc.edu.cn}
		\address[mymainaddress1]{Department of Statistics, Zhejiang University City College, Hangzhou,
			310015, China}
		\begin{abstract}
			In this paper, we study the sharp oracle bounds for Slope and Lasso and generalize the results in \cite{bellec2018} to allow the case that the parameter vector is not exactly sparse and obtain the optimal bounds for $\ell_q$ estimation errors with $1\leq q\leq \infty$ by using some extended Restricted Eigenvalue type conditions.
		\end{abstract}
		
		
	\end{frontmatter}

\section{Introduction}

In high-dimensional linear regression \cite{buhlmann2011}, we aim to estimate the $s$-sparse parameter vector $\beta^{*}\in\mathbb{R}^p$ (i.e., $\lVert\beta^{*}\rVert_0\leq s$) from few noisy observations $y=X\beta^{*}+\xi$, where $y\in\mathbb{R}^n(n\ll p)$, $X\in\mathbb{R}^{n\times p}$ is the design matrix and $\xi\in\mathbb{R}^n\sim N(0,\sigma^2 I_{n\times n})$ is the Gaussian noise. The best prediction and $\ell_2$ estimation rate achieved by the well-known polynomial time algorithms such as the Lasso \cite{tibshirani1996} and the Dantzig selector \cite{bickel2009} is $(s/n)\log (p)$, while the minimax rate for the problem is $(s/n)\log(p/s)$ \cite{raskutti2011}. Until recently, \cite{su2016} has shown that for i.i.d. standard normal regressor this optimal rate can be asymptotically achieved by a computable in polynomial time method called Sorted L-One Penalized Estimator (Slope) \cite{bogdan2015}. For the high-dimensional linear regression problem, Slope is defined as a solution of the minimization problem: \begin{align}
\mathop{\arg\min}\limits_{\beta\in\mathbb{R}^p} \frac{1}{n}\lVert y-X\beta\rVert_2^2+2\lVert\beta\rVert_{*}, \label{slope_definition}
\end{align}
where the norm $\lVert\beta\rVert_{*}=\sum\limits_{j=1}^p \lambda_j\beta_j^{\sharp}$, the tuning parameters satisfy $\lambda_1\geq \lambda_2\geq \cdots\lambda_p\geq 0$, $(\beta_1^{\sharp},\beta_2^{\sharp},\cdots,\beta_p^{\sharp})$ is a non-increasing rearrangement of $(|\beta_1|,|\beta_2|,\cdots,|\beta_p|)$. \cite{lecue2018} extends the results in \cite{su2016} to the subgaussian designs and obtains a non-asymptotic bound with this optimal rate by using the small ball method. Nevertheless, the strong assumption that the covariance matrix of the design is proportional to the identity matrix (namely the design is isotropic) is required in both of these works. In order to relax this restriction, \cite{bellec2018} shows that under the Restricted Eigenvalue (RE) type conditions a Lasso estimator with adaptively chosen tuning parameters and a Slope estimator can achieve the non-asymptotic minimax prediction and $\ell_2$ estimation rate $(s/n)\log(p/s)$ both in probability and in expectation. The RE type conditions imposed on the design are satisfied for a large class of random matrices with independent and possibly anisotropically distributed rows with high probability. The Slope method has attracted more and more research interest. Among them, \cite{derumigny2018} discussed the Squared-Root Slope method. \cite{brzyski2019} extended to the Group Slope for handling the group sparse regression problem, while \cite{feng2019} proposed the Sorted Concave Penalized method through combining the Slope and Concave Penalized method.

However, we notice that all the existing literature on the Slope only discussed the case that the parameter vector to be estimated is exactly sparse and the bounds are merely obtained for the $\ell_q$ estimation errors with $1\leq q\leq 2$. In this paper, we aim to generalize the results in \cite{bellec2018} to allow the case that the parameter vector is not exactly sparse and derive the optimal bounds for $\ell_q$ estimation errors with $1\leq q\leq \infty$ via the extended RE type conditions. 

The paper is organized as follows. In Section 2, we study the estimation performance of the Lasso estimator. In Section 3, we present the results for the Slope estimator. In Section 4, we discuss the optimality of the derived upper bounds. In Section 5, we investigate the proposed RE conditions. In Section 6, we give the conclusion.

\section{Lasso estimator}
In this section, we denote the Lasso estimator by $\hat{\beta}\in\mathop{\arg\min}\limits_{\beta\in\mathbb{R}^p} \frac{1}{n}\lVert y-X\beta\rVert_2^2+2\lambda\lVert\beta\rVert_{1}$, where $\lambda>0$ is a tuning parameter. We use some notations in \cite{bellec2018} in order to facilitate the presentation of the results. For a given $\delta_0\in(0,1)$ and for any $u=(u_1,u_2,\cdots,u_p)\in\mathbb{R}^p$, we set \begin{align}
H(u):=(4+\sqrt{2})\sum\limits_{j=1}^p u_j^{\sharp}\sigma\sqrt{\frac{\log(2p/j)}{n}},\,\, G(u):=(4+\sqrt{2})\sigma\sqrt{\frac{\log(1/\delta_0)}{n}}\lVert Xu\rVert_2/\sqrt{n}, \label{Notation}
\end{align}
where $(u_1^\sharp,u_2^\sharp,\cdots,u_p^\sharp)$ is a non-increasing rearrangement of $(|u_1|,|u_2|,\cdots,|u_p|)$. For a vector $\beta\in\mathbb{R}^p$ and a set $S\subset\{1,2,\cdots,p\}$, we denote by $\beta_S$ the vector which
coincides with $\beta$ on the indices in $S$ and is extended to zero outside $S$.

First, as an extension of LEMMA A.1 in \cite{bellec2018}, we have the following key lemma that allows the parameter vector to be not exactly sparse. 

\begin{lemma}
	Let $s\in\{1,\cdots,p\}$, $\tau\in[0,1]$, $\lambda_1\geq \lambda_2\geq\cdots\geq\lambda_p\geq 0$ and the error of the best $s$-term approximation of $\beta$ with respect to $\lVert\cdot\rVert_{*}$ is defined by $\sigma_s(\beta)_{*}=\min\limits_{\lVert z\rVert_0\leq s}\lVert \beta-z\rVert_{*}$. For any two $\beta, \hat{\beta}\in\mathbb{R}^p$ and $1<q\leq \infty$, we have \begin{align}
	\tau\lVert u\rVert_{*}+\lVert \beta\rVert_{*}-\lVert \hat{\beta}\rVert_{*}&\leq (1+\tau)\sum\limits_{j=1}^s \lambda_j u_j^{\sharp}-(1-\tau)\sum\limits_{j=s+1}^p \lambda_j u_j^{\sharp}+2\sigma_s(\beta)_{*}  \nonumber \\
	&\leq (1+\tau)\left(\sum\limits_{j=1}^s \lambda_j^{\frac{q}{q-1}}\right)^{1-1/q}\lVert u\rVert_q-(1-\tau)\sum\limits_{j=s+1}^p \lambda_j u_j^{\sharp}+2\sigma_s(\beta)_{*}, \label{slope}
	\end{align}
	where $u=\hat{\beta}-\beta=(u_1,u_2,\cdots,u_p)^T$ and $(u_{1}^{\sharp},u_2^\sharp,\cdots,u_p^{\sharp})$ is a non-increasing rearrangement of $(|u_1|,|u_2|,\cdots,|u_p|)$. In particular, if $\lambda_1=\lambda_2=\cdots=\lambda_p=\lambda$ for some $\lambda>0$, then $\lVert \cdot\rVert_{*}=\lambda\lVert\cdot\rVert_1$ and we have \begin{align}
	\tau\lambda\lVert u\rVert_1+\lambda\lVert\beta\rVert_1-\lambda\lVert\hat{\beta}\rVert_1 &\leq (1+\tau)\lambda\sum\limits_{j=1}^s u_j^{\sharp}-(1-\tau)\lambda\sum\limits_{j=s+1}^p  u_j^{\sharp}+2\lambda\sigma_s(\beta)_1 \nonumber \\
	&\leq (1+\tau)\lambda s^{1-1/q}\lVert u\rVert_q-(1-\tau)\lambda\sum\limits_{j=s+1}^p  u_j^{\sharp}+2\lambda\sigma_s(\beta)_1.
	\end{align}
\end{lemma}

\noindent
{\bf Proof.} Let $S$ be the index set of $s$ largest absolute entries of $\beta$ so that $\sigma_s(\beta)_{*}=\lVert \beta_{S^c}\rVert_{*}$. Assume $\phi$ is any permutation of $\{1,\cdots,p\}$ such that \begin{align*}
\lVert \beta_S\rVert_{*}=\sum\limits_{j=1}^s \lambda_j |\beta_{\phi(j)}|\,\,\text{and $|u_{\phi(s+1)}|\geq |u_{\phi(s+2)}|\geq\cdots\geq|u_{\phi(p)}|$.}
\end{align*}
Then by using the definition of $\lVert\cdot\rVert_{*}$ and the Triangle inequality, we have \begin{align*}
\lVert \beta\rVert_{*}-\lVert \hat{\beta}\rVert_{*}&\leq \lVert \beta_S\rVert_{*}+\lVert \beta_{S^c}\rVert_{*}-\lVert \hat{\beta}\rVert_{*} \\
&\leq \sum\limits_{j=1}^s \lambda_j \left(|\beta_{\phi(j)}|-|\hat{\beta}_{\phi(j)}|\right)-\sum\limits_{j=s+1}^p \lambda_j|\hat{\beta}_{\phi(j)}|+\lVert\beta_{S^c}\rVert_{*} \\
&\leq \sum\limits_{j=1}^s \lambda_j|u_{\phi(j)}|-\sum\limits_{j=s+1}^p \lambda_j|\hat{\beta}_{\phi(j)}|+\lVert\beta_{S^c}\rVert_{*} \\
&\leq \sum\limits_{j=1}^s \lambda_j|u_{\phi(j)}|-\sum\limits_{j=s+1}^p \lambda_j|u_{\phi(j)}|+\sum\limits_{j=s+1}^p \lambda_j|\beta_{\phi(j)}|+\lVert\beta_{S^c}\rVert_{*} \\
&\leq \sum\limits_{j=1}^s \lambda_j|u_{\phi(j)}|-\sum\limits_{j=s+1}^p \lambda_j|u_{\phi(j)}|+2\lVert\beta_{S^c}\rVert_{*} \\
&\leq \sum\limits_{j=1}^s \lambda_j u_j^{\sharp}-\sum\limits_{j=s+1}^p \lambda_j u_j^{\sharp}+2\sigma_s(\beta)_{*},
\end{align*}
since $\sum\limits_{j=s+1}^p \lambda_j|\beta_{\phi(j)}|\leq \lVert\beta_{S^c}\rVert_{*}$, $\sum\limits_{j=1}^s \lambda_j|u_{\phi(j)}|\leq \sum\limits_{j=1}^s \lambda_j u_j^{\sharp}$ and $\sum\limits_{j=s+1}^p \lambda_j u_j^{\sharp}\leq\sum\limits_{j=s+1}^p \lambda_j|u_{\phi(j)}|$. Finally, (\ref{slope}) follows from $\sum\limits_{j=1}^s \lambda_j u_j^{\sharp}\leq \left(\sum\limits_{j=1}^s \lambda_j^{\frac{q}{q-1}}\right)^{1-1/q}\lVert u\rVert_q$ by the H\"{o}lder's inequality.\\

To be self-contained, here we list two technical lemmas given in \cite{bellec2018}, which will be used for establishing the main results in what follows.

\begin{lemma}(Lemma A.2 in \cite{bellec2018})
	Let $h:\mathbb{R}^p\rightarrow \mathbb{R}$ be a convex function, let $f,\xi\in\mathbb{R}^n$, $y=f+\xi$ and let $X$ be any $n\times p$ matrix. If $\hat{\beta}$ is a solution of the minimization problem $\min_{\beta\in\mathbb{R}^p}\{\frac{1}{n}\lVert X\beta-y\rVert_2^2+h(\beta)\}$, then $\hat{\beta}$ satisfies for all $\beta\in\mathbb{R}^p$ \begin{align}
	\frac{1}{n}\lVert X\hat{\beta}-f\rVert_2^2-\frac{1}{n}\lVert X\beta-f\rVert_2^2\leq \frac{2}{n}\xi^T X(\hat{\beta}-\beta)+h(\beta)-h(\hat{\beta})-\frac{1}{n}\lVert X(\hat{\beta}-\beta)\rVert_2^2. \label{KKT}
	\end{align}
\end{lemma}
\begin{lemma} (Theorem 4.1 in \cite{bellec2018})
	Let $\delta_0\in (0,1)$ and let $X\in\mathbb{R}^{n\times p}$ be a matrix such that $\max_{j=1,\cdots,p}\lVert Xe_j\rVert_2/\sqrt{n}\leq 1$ where $(e_1,e_2,\cdots,e_p)$ is the canonical basis in $\mathbb{R}^p$. Let $H(\cdot)$ and $G(\cdot)$ be defined in (\ref{Notation}). If $\xi\sim N(0,\sigma^2 I_{n\times n})$, then the random event \begin{align}
	\left\{\frac{1}{n}\xi^{T} X u \leq \max (H(u),G(u)), \forall\,u\in\mathbb{R}^p\right\} \label{event}
	\end{align}
	is of probability of at least $1-\delta_0/2$.
\end{lemma}

In order to obtain the bounds for $\ell_q$ estimation errors with any $q\in[1,\infty]$, we generalize the Strong Restricted Eigenvalue (SRE) condition in \cite{bellec2018} to involve the variable $q$ lying in $[2,\infty]$. When $q=2$, the following $SRE(q,s,c_0)$ condition reduces to the $SRE(s,c_0)$ condition used in \cite{bellec2018}. In fact, we notice that the SRE condition is exactly equivalent to the constrained minimal singular values (CMSV) condition proposed in \cite{tang2011}. Thus, it's natural to propose the following general $SRE(q,s,c_0)$ condition, which is  essentially the same as the $q$-ratio CMSV condition studied in \cite{zhou2019a,zhou2019b} for the sparse signal recovery. For any $s\in[1,p]$ and $q\in (1,\infty]$, the $q$-ratio CMSV of a design matrix $X\in\mathbb{R}^{n\times p}$ is defined as $\rho_{q,s}(X)=\min_{\{\delta\neq 0, s_q(\delta)\leq s\}}\frac{\lVert X\delta\rVert_2}{\lVert \delta\rVert_q}$ where $s_q(\delta)=(\frac{\lVert\delta\rVert_1}{\lVert\delta\rVert_q})^{\frac{q}{q-1}}$ is the $q$-ratio sparsity level of non-zero $\delta\in\mathbb{R}^p$.

\begin{definition}($SRE(q,s,c_0)$ condition)
Let $c_0>0$ and $s\in\{1,2,\cdots,p\}$. For any $2\leq q\leq \infty$, the design matrix $X$ satisfies $\lVert Xe_j\rVert_2/\sqrt{n}\leq 1$ for all $j=1,\cdots,p$ and \begin{align}
\theta_q(s,c_0):=\min\limits_{\delta\neq 0, \delta\in\mathcal{C}_{SRE(q,s,c_0)}} \frac{\lVert X\delta\rVert_2}{\sqrt{n}\lVert \delta\rVert_q}>0,
\end{align}
where $\mathcal{C}_{SRE(q,s,c_0)}:=\{\delta\in\mathbb{R}^p: \lVert \delta\rVert_1\leq (1+c_0) s^{1-1/q}\lVert \delta\rVert_q\}$ is a cone in $\mathbb{R}^p$. 
\end{definition}

Under the new proposed SRE condition, we now present our main result for the Lasso estimator. Let $\gamma\in(0,1)$ be a constant. For any tuning parameter $\lambda>0$, we set \begin{align}
\delta(\lambda):=\exp\left(-\left(\frac{\gamma\lambda\sqrt{n}}{(4+\sqrt{2})\sigma}\right)^2\right)\,\,\text{so that $\lambda=\frac{(4+\sqrt{2})\sigma}{\gamma}\sqrt{\frac{\log(1/\delta(\lambda))}{n}}$.}
\end{align}
For a given $s\in\{1,\cdots,p\}$, the following theorem holds under the condition \begin{align}
\lambda\geq \frac{(4+\sqrt{2})\sigma}{\gamma}\sqrt{\frac{\log(2ep/s)}{n}}\,\,\text{or equivalently $\delta(\lambda)\leq s/(2ep)$}. \label{lambda} 
\end{align}

\begin{theorem}
	Let $s\in\{1,\cdots,p\}$, $\gamma\in(0,1)$ and $\tau\in[0, 1-\gamma)$. For any $2\leq q\leq \infty$, assume that the $SRE(q,s,c_0)$ condition holds with $c_0=c_0(\gamma,\tau)=\frac{1+\gamma+\tau}{1-\gamma-\tau}$. Let $\lambda$ be a tuning parameter such that (\ref{lambda}) holds. Let $\delta_0\in(0,1)$. Then, on the event (\ref{event}), the Lasso estimator $\hat{\beta}$ with tuning parameter $\lambda$ satisfies \begin{align}
	2\tau\lambda\lVert\hat{\beta}-\beta\rVert_1+\frac{1}{n}\lVert X\hat{\beta}-f\rVert_2^2\leq \frac{1}{n}\lVert X\beta-f\rVert_2^2+C_{\gamma,\tau}(q,s,\lambda,\delta_0)\lambda^2s+4\lambda\sigma_s(\beta)_1, \label{result1}
	\end{align}
	for all $\beta\in\mathbb{R}^p$ and all $f\in\mathbb{R}^n$, where \begin{align*}
	C_{\gamma,\tau}(q,s,\lambda,\delta_0):=(1+\gamma+\tau)^2\left(\frac{\log(1/\delta_0)}{s\log(1/\delta(\lambda))}\vee \frac{s^{1-2/q}}{\theta_q^2(s,c_0(\gamma,\tau))}\right).
	\end{align*}
	Furthermore, if $f=X\beta^{*}$ for some $\beta^{*}$, then on the event (\ref{event}), we have \begin{align}
	\lVert \hat{\beta}-\beta^{*}\rVert_1&\leq \frac{C_{\gamma,\tau}(q,s,\lambda,\delta_0)}{2\tau}\lambda s+\frac{2}{\tau}\sigma_s(\beta^{*})_1,  \label{result2} \\ 
	\lVert \hat{\beta}-\beta^{*}\rVert_q&\leq \frac{C_{\gamma,0}(q,s,\lambda,\delta_0)}{1+\gamma}\lambda s^{1/q},\quad (\text{when $\beta^{*}$ is $s$-sparse}) \label{result3} \\
	\lVert \hat{\beta}-\beta^{*}\rVert_q&\leq \left(\frac{2}{1+\gamma}\vee \frac{1-\gamma-\tau}{4\tau}\right)C_{\gamma,\tau}(q,s,\lambda,\delta_0)\lambda s^{1/q}+\left(\frac{1}{1+\gamma}\vee \frac{1-\gamma-\tau}{\tau}\right)s^{1/q-1}\sigma_s(\beta^{*})_1. \label{result4}
	\end{align}
\end{theorem}

\medskip
\noindent
{\bf Proof.} By using (\ref{KKT}) with $h(\cdot)=2\lambda\lVert \cdot\rVert_1$, we get that for all $\beta\in\mathbb{R}^p$ and all $f\in\mathbb{R}^n$ \begin{align}
2\tau\lambda\lVert \hat{\beta}-\beta\rVert_1+\frac{1}{n}\lVert X\hat{\beta}-f\rVert_2^2\leq \frac{1}{n}\lVert X\beta-f\rVert_2^2-\frac{1}{n}\lVert X(\hat{\beta}-\beta)\rVert_2^2+\triangle^{*}, \label{inequality}
\end{align}
where \begin{align*}
\triangle^{*}:=2\tau\lambda\lVert\hat{\beta}-\beta\rVert_1+\frac{2}{n}\xi^T X(\hat{\beta}-\beta)+2\lambda\lVert\beta\rVert_1-2\lambda\lVert \hat{\beta}\rVert_1.
\end{align*}
We let $u=\hat{\beta}-\beta$ and define \begin{align}
\tilde{H}(u):=\frac{(4+\sqrt{2})\sigma}{\sqrt{n}}\left(\lVert u\rVert_q \left[\sum\limits_{j=1}^s\left(\sqrt{\log (2p/j)}\right)^{\frac{q}{q-1}}\right]^{1-1/q}+\sum\limits_{j=s+1}^p u_j^{\sharp}\sqrt{\log(2p/j)}\right).
\end{align}
For any $2\leq q\leq \infty$, using the H\"{o}lder's inequality and the fact that $\sum\limits_{j=1}^s \log(2p/j)\leq s\log(2ep/s)$, it's easy to verify that \begin{align}
H(u)=\frac{(4+\sqrt{2})\sigma}{\sqrt{n}}\sum\limits_{j=1}^p u_j^\sharp\sqrt{\log(2p/j)}\leq \tilde{H}(u)\leq \lambda\gamma\left(s^{1-1/q}\lVert u\rVert_q+\sum\limits_{j=s+1}^p u_j^{\sharp}\right)=: F(u), \forall\, u\in\mathbb{R}^p.
\end{align}

On the event (\ref{event}), by using Lemma 1 we have \begin{align*}
\triangle^{*}&\leq 2\lambda(\tau \lVert \hat{\beta}-\beta\rVert_1+\lVert \beta\rVert_1-\lVert\hat{\beta}\rVert_1)+2\max(H(u),G(u)) \\
&\leq 2\lambda\left((1+\tau)s^{1-1/q}\lVert u\rVert_q-(1-\tau)\sum\limits_{j=s+1}^p u_j^{\sharp}+2\sigma_s(\beta)_1\right)+2\max(F(u),G(u)).
\end{align*}

By the definition of $\delta(\lambda)$, i.e., $\lambda=\frac{(4+\sqrt{2})\sigma}{\gamma}\sqrt{\frac{\log(1/\delta(\lambda))}{n}}$, we have $$
G(u)=\lambda\sqrt{s}\gamma\sqrt{\frac{\log(1/\delta_0)}{s\log(1/\delta(\lambda))}}\lVert Xu\rVert_2/\sqrt{n}.
$$
As the arguments in Proof of Theorem 4.2 in \cite{bellec2018}, we consider the following two cases: \\

(i) Case $G(u)>F(u)$: Then it holds that \begin{align}
\lVert u\rVert_q\leq s^{1/q-1/2}\sqrt{\frac{\log(1/\delta_0)}{s\log(1/\delta(\lambda))}}\lVert Xu\rVert_2/\sqrt{n}. \label{ratio1}
\end{align}
Thus, \begin{align}
\triangle^{*}&\leq 2\lambda(1+\tau)s^{1-1/q}\lVert u\rVert_q+2G(u)+4\lambda\sigma_s(\beta)_1 \nonumber \\
&\leq 2\lambda(1+\tau)\sqrt{s}\sqrt{\frac{\log(1/\delta_0)}{s\log(1/\delta(\lambda))}}\lVert Xu\rVert_2/\sqrt{n}+2\lambda\sqrt{s}\gamma\sqrt{\frac{\log(1/\delta_0)}{s\log(1/\delta(\lambda))}}\lVert Xu\rVert_2/\sqrt{n}+4\lambda\sigma_s(\beta)_1 \nonumber  \\
&\leq 2\lambda\sqrt{s}(1+\tau+\gamma)\sqrt{\frac{\log(1/\delta_0)}{s\log(1/\delta(\lambda))}}\lVert Xu\rVert_2/\sqrt{n}+4\lambda\sigma_s(\beta)_1 \nonumber \\
&\leq \lambda^2s(1+\tau+\gamma)^2\frac{\log(1/\delta_0)}{s\log(1/\delta(\lambda))}+\frac{1}{n}\lVert Xu\rVert_2^2+4\lambda\sigma_s(\beta)_1.  \label{triangle1} 
\end{align}

(ii) Case $G(u)\leq F(u)$:  In this case it follows that \begin{align*}
\triangle^{*}\leq \underbrace{2\lambda\left((1+\gamma+\tau)s^{1-1/q}\lVert u\rVert_q-(1-\gamma-\tau)\sum\limits_{j=s+1}^p u_j^{\sharp}\right)}_{\triangle}+4\lambda\sigma_s(\beta)_1.
\end{align*}
If $\triangle>0$, then $\lVert u\rVert_1=\sum\limits_{j=1}^s u_j^\sharp+\sum\limits_{j=s+1}^p u_j^\sharp\leq (1+c_0(\gamma,\tau))s^{1-1/q}\lVert u\rVert_q$, i.e., $u$ belongs to the cone $\mathcal{C}_{SRE}(q,s,c_0)$. The $SRE(q,s,c_0)$ condition yields $\lVert u\rVert_q\leq \frac{\lVert Xu\rVert_2}{\theta_q(s,c_0)\sqrt{n}}$. Therefore, \begin{align}
\triangle^{*}&\leq \triangle+4\lambda\sigma_s(\beta)_1\leq \frac{2\lambda(1+\gamma+\tau)s^{1-1/q}\lVert Xu\rVert_2/\sqrt{n}}{\theta_q(s,c_0)}+4\lambda\sigma_s(\beta)_1 \nonumber\\
&\leq \left(\frac{\lambda(1+\gamma+\tau)s^{1-1/q}}{\theta_q(s,c_0)}\right)^2+\frac{1}{n}\lVert Xu\rVert_2^2+4\lambda\sigma_s(\beta)_1. \label{triangle2}
\end{align}
If $\triangle\leq 0$, then (\ref{triangle2}) holds trivially. 

The proof of (\ref{result1}) is completed by combining (\ref{triangle1}) and (\ref{triangle2}) with (\ref{inequality}). 

Now let $f=X\beta^{*}$ and $\beta=\beta^{*}$ in (\ref{result1}), then \begin{align*}
2\tau\lambda\lVert\hat{\beta}-\beta^{*}\rVert_1\leq2\tau\lambda\lVert\hat{\beta}-\beta^{*}\rVert_1+\frac{1}{n}\lVert X(\hat{\beta}-\beta^{*})\rVert_2^2 \leq C_{\gamma,\tau}(q,s,\lambda,\delta_0)\lambda^2 s+4\lambda\sigma_s(\beta^{*})_1,
\end{align*}
which immediately implies (\ref{result2}). 
 
Next, to prove (\ref{result3}), when the true parameter vector $\beta^{*}$ is exactly $s$-sparse (i.e., $\sigma_s(\beta^{*})_1=0$), we take $\tau=0$ and consider the Cases (i) and (ii) as above with $u=\hat{\beta}-\beta^{*}$.  \begin{itemize}
 	\item If $G(u)>F(u)$, then (\ref{inequality}) and (\ref{triangle1}) implies that $$
 	\frac{1}{n}\lVert X(\hat{\beta}-\beta^{*})\rVert_2^2\leq \lambda^2(1+\gamma)^2\frac{\log(1/\delta_0)}{\log(1/\delta(\lambda))}. 
 	$$
By this and (\ref{ratio1}), we get \begin{align}
\lVert \hat{\beta}-\beta^{*}\rVert_q&\leq s^{1/q-1/2}\sqrt{\frac{\log(1/\delta_0)}{s\log(1/\delta(\lambda))}}\cdot\sqrt{\lambda^2(1+\gamma)^2\frac{\log(1/\delta_0)}{\log(1/\delta(\lambda))}} \nonumber \\
&\leq s^{1/q-1/2}\sqrt{\frac{C_{\gamma,0}(q,s,\lambda,\delta_0)}{(1+\gamma)^2}}\cdot
\sqrt{\lambda^2(1+\gamma)^2\frac{sC_{\gamma,0}(q,s,\lambda,\delta_0)}{(1+\gamma)^2}}  \nonumber \\
&\leq \frac{C_{\gamma,0}(q,s,\lambda,\delta_0)}{1+\gamma}\lambda s^{1/q}. \label{bound1}
 \end{align}
 \item If $G(u)\leq F(u)$, we have $\triangle\geq \triangle^{*}\geq 0$, i.e., $u=\hat{\beta}-\beta^{*}\in \mathcal{C}_{SRE(q,s,\frac{1+\gamma}{1-\gamma})}$ such that \begin{align*}
 \frac{1}{n}\lVert Xu\rVert_2^2\leq \left(\frac{(1+\gamma)\lambda s^{1-1/q}}{\theta_q(s,\frac{1+\gamma}{1-\gamma})}\right)^2.
 \end{align*} 
Thus we can apply the $SRE(q,s,\frac{1+\gamma}{1-\gamma})$ condition, which yields that \begin{align}
 \lVert \hat{\beta}-\beta^{*}\rVert_q\leq \frac{\lVert Xu\rVert_2/\sqrt{n}}{\theta_q(s,\frac{1+\gamma}{1-\gamma})}&\leq \frac{1}{\theta_q(s,\frac{1+\gamma}{1-\gamma})}\sqrt{\left(\frac{(1+\gamma)\lambda s^{1-1/q}}{\theta_q(s,\frac{1+\gamma}{1-\gamma})}\right)^2} \nonumber \\
& \leq \frac{\sqrt{C_{\gamma,0}(q,s,\lambda,\delta_0)}}{1+\gamma}s^{1/q-1/2}\sqrt{\lambda^2 sC_{\gamma,0}(q,s,\lambda,\delta_0)} \nonumber \\
& \leq \frac{C_{\gamma,0}(q,s,\lambda,\delta_0)}{1+\gamma}\lambda s^{1/q}. \label{bound2}
 \end{align}
 \end{itemize}
Putting together (\ref{bound1}) and (\ref{bound2}) yields $
 \lVert \hat{\beta}-\beta^{*}\rVert_q\leq \frac{C_{\gamma,0}(q,s,\lambda,\delta_0)}{1+\gamma}\lambda s^{1/q}$,
which completes the proof of (\ref{result3}).

As for the result (\ref{result4}) that allows the true parameter vector $\beta^{*}$ to be not exactly sparse, we let $f=X\beta^{*}$, $\beta=\beta^{*}$ and $u=\hat{\beta}-\beta^{*}$ in (\ref{result1}), then \begin{align}
2\tau\lambda\lVert u\rVert_1+\frac{1}{n}\lVert Xu\rVert_2^2\leq \lambda^2 sC_{\gamma,\tau}(q,s,\lambda,\delta_0)+4\lambda\sigma_s(\beta^{*})_1.
\end{align}
If $u\in\mathcal{C}_{SRE}(q,s,c_0=c_0(\gamma,\tau)=\frac{1+\gamma+\tau}{1-\gamma-\tau})$, then the $SRE(q,s,c_0)$ condition implies that \begin{align}
\lVert u \rVert_q\leq \frac{\lVert Xu\rVert_2/\sqrt{n}}{\theta_q(s,c_0(\gamma,\tau))}&\leq \frac{1}{\theta_q(s,c_0(\gamma,\tau))}\sqrt{\lambda^2 sC_{\gamma,\tau}(q,s,\lambda,\delta_0)+4\lambda\sigma_s(\beta^*)_1} \nonumber \\
& \leq \frac{\sqrt{C_{\gamma,\tau}(q,s,\lambda,\delta_0)}}{1+\gamma}s^{1/q-1/2}\sqrt{\lambda^2 sC_{\gamma,\tau}(q,s,\lambda,\delta_0)+4\lambda\sigma_s(\beta^*)_1} \nonumber \\
& \leq \frac{C_{\gamma,\tau}(q,s,\lambda,\delta_0)}{1+\gamma}\lambda s^{1/q}+\frac{1}{1+\gamma}\left(2\sqrt{\lambda C_{\gamma,\tau}(q,s,\lambda,\delta_0) }s^{\frac{1}{2q}}\cdot s^{\frac{1}{2}(\frac{1}{q}-1)}\sqrt{\sigma_s(\beta^*)_1}\right) \nonumber \\
& \leq \frac{C_{\gamma,\tau}(q,s,\lambda,\delta_0)}{1+\gamma}\lambda s^{1/q}+\frac{1}{1+\gamma}\left(\lambda C_{\gamma,\tau}(q,s,\lambda,\delta_0) s^{1/q}+s^{1/q-1}\sigma_s(\beta^*)_1\right) \nonumber \\
& \leq \frac{2}{1+\gamma}\lambda C_{\gamma,\tau}(q,s,\lambda,\delta_0) s^{1/q}+\frac{1}{1+\gamma}s^{1/q-1}\sigma_s(\beta^*)_1. \label{compressible1}
\end{align}
Otherwise, namely $u\notin\mathcal{C}_{SRE}(q,s,c_0=c_0(\gamma,\tau)=\frac{1+\gamma+\tau}{1-\gamma-\tau})$, by using (\ref{result2}) we have \begin{align}
\lVert u\rVert_q<\frac{1}{1+c_0}s^{1/q-1}\lVert u\rVert_1&<\frac{1-\gamma-\tau}{2}s^{1/q-1}\left(\frac{C_{\gamma,\tau}(q,s,\lambda,\delta_0)}{2\tau}\lambda s+\frac{2}{\tau}\sigma_s(\beta^{*})_1\right) \nonumber \\
&<\frac{1-\gamma-\tau}{4\tau}\lambda C_{\gamma,\tau}(q,s,\lambda,\delta_0) s^{1/q}+\frac{1-\gamma-\tau}{\tau}s^{1/q-1}\sigma_s(\beta^*)_1.\label{compressible2}
\end{align}
Then, (\ref{result4}) follows immediately from (\ref{compressible1}) and (\ref{compressible2}). \\

By letting $q=2$ in Theorem 1, and adopting the norm interpolation inequality $\lVert \hat{\beta}-\beta^{*}\rVert_q\leq \lVert \hat{\beta}-\beta^{*}\rVert_1^{2/q-1}\lVert \hat{\beta}-\beta^{*}\rVert_2^{2-2/q}$ for any $1\leq q\leq 2$, we can obtain the following corollary immediately, which is exactly the Theorem 4.2 in \cite{bellec2018}. 

\begin{cor}
 Under the same conditions of Theorem 1 and for any $1\leq q\leq 2$, then on the event (\ref{event}), if $\beta^*$ is $s$-sparse, we have \begin{align}
 \lVert \hat{\beta}-\beta^{*}\rVert_q\leq \left(\frac{C_{\gamma,\tau}(2,s,\lambda,\delta_0)}{2\tau}\right)^{2/q-1}\left(\frac{C_{\gamma,0}(2,s,\lambda,\delta_0)}{1+\gamma}\right)^{2-2/q}\lambda s^{1/q}.
 \end{align}
\end{cor}

In addition, if we use the particular choice of constants with $\gamma=\frac{1}{2}$, $\tau=\frac{1}{4}$, $\delta_0=\delta_0^{*}=(\delta(\lambda))^{\frac{s}{\theta_{2}^2(s,7)}}$ for Corollary 1 while  $\delta_0=\delta_0^{*}=(\delta(\lambda))^{\frac{s^{2-2/q}}{\theta_{q}^2(s,7)}}\leq \left(\frac{s}{2ep}\right)^{\frac{s^{2-2/q}}{\theta_{q}^2(s,7)}}$ for Theorem 2, then we are able to derive the Corollary 4.3 in \cite{bellec2018} and the corollary as follows.

\begin{cor}
	Let $s\in\{1,\cdots,p\}$. For any $2\leq q\leq \infty$, assume that the $SRE(q,s,7)$ condition holds. Let $\lambda$ be a tuning parameter such that (\ref{lambda}) holds with $\gamma=1/2$. Then, with probability at least $1-\frac{1}{2}\left(\frac{s}{2ep}\right)^{\frac{s^{2-2/q}}{\theta_{q}^2(s,7)}}$, the Lasso estimator $\hat{\beta}$ with tuning parameter $\lambda$ satisfies \begin{align*}
	\frac{\lambda}{2}\lVert\hat{\beta}-\beta\rVert_1+\frac{1}{n}\lVert X\hat{\beta}-f\rVert_2^2\leq \frac{1}{n}\lVert X\beta-f\rVert_2^2+\frac{49\lambda^2s^{2-2/q}}{16\theta_{q}^2(s,7)}+4\lambda\sigma_s(\beta)_1, 
	\end{align*}
	for all $\beta\in\mathbb{R}^p$ and all $f\in\mathbb{R}^n$. Furthermore, if $f=X\beta^{*}$ for some $\beta^{*}$, then with probability at least $1-\frac{1}{2}\left(\frac{s}{2ep}\right)^{\frac{s^{2-2/q}}{\theta_{q}^2(s,7)}}$, we have \begin{align*}
	\lVert \hat{\beta}-\beta^{*}\rVert_1&\leq \frac{49\lambda s^{2-2/q}}{8\theta_q^2(s,7)}+8\sigma_s(\beta^{*})_1,   \\ 
	\lVert \hat{\beta}-\beta^{*}\rVert_q&\leq \frac{49\lambda s^{1-1/q}}{24\theta_q^2(s,7)},\quad (\text{when $\beta^{*}$ is $s$-sparse})  \\
	\lVert \hat{\beta}-\beta^{*}\rVert_q&\leq  \frac{49\lambda s^{1-1/q}}{12\theta_q^2(s,7)}+s^{1/q-1}\sigma_s(\beta^{*})_1. 
	\end{align*}
\end{cor}
	
\section{Slope estimator}

In this section, we study the Slope estimator $\hat{\beta}$ given as a solution of the minimization problem (\ref{slope_definition}) with weights \begin{align}
\lambda_j=A\sigma \sqrt{\frac{\log(2p/j)}{n}} \label{weights}
\end{align} for any constant $A>4+\sqrt{2}$. Hereafter, we let $\Lambda_q(s)=\left(\sum_{j=1}^s \lambda_j^\frac{q}{q-1}\right)^{1-1/q}$. We use the following condition on the design matrix $X$ called $WRE(q,s,c_0)$ condition with $2\leq q\leq \infty$, which acts as an extension of the Weighted Restricted Eigenvalue (WRE) condition used in \cite{bellec2018}.

\begin{definition}($WRE(q,s,c_0)$ condition)
	Let $c_0>0$ and $s\in\{1,2,\cdots,p\}$. For any $2\leq q\leq \infty$, the design matrix $X$ satisfies $\lVert Xe_j\rVert_2/\sqrt{n}\leq 1$ for all $j=1,\cdots,p$ and \begin{align}
	\nu_q(s,c_0):=\min\limits_{\delta\neq 0, \delta\in\mathcal{C}_{WRE(q,s,c_0)}} \frac{\lVert X\delta\rVert_2}{\sqrt{n}\lVert \delta\rVert_q}>0,
	\end{align}
	where $\mathcal{C}_{WRE(q,s,c_0)}:=\{\delta\in\mathbb{R}^p: \lVert \delta\rVert_{*}\leq (1+c_0) \lVert \delta\rVert_q\Lambda_q(s)\}$ is a cone in $\mathbb{R}^p$. 
\end{definition}

\medskip
This $WRE(q,s,c_0)$ condition is defined for any weights $\lambda_1\geq \cdots\geq \lambda_p\geq 0$, while we only focus on the weights given in (\ref{weights}) here. In this case, the cone is equivalently defined as\[ \mathcal{C}_{WRE(q,s,c_0)}:=\left\{\delta\in\mathbb{R}^p: \sum\limits_{j=1}^p \delta_j^\sharp \sqrt{\log(2p/j)}\leq (1+c_0) \lVert \delta\rVert_q \left[ \sum\limits_{j=1}^s (\sqrt{\log(2p/j)})^{\frac{q}{q-1}}\right]^{1-1/q}\right\}.
\]
Thus, if $\delta\in \mathcal{C}_{SRE(q,s,c_0)}$ (i.e., $\lVert \delta\rVert_1\leq (1+c_0)s^{1-1/q}\lVert \delta\rVert_q$), we have $\sum\limits_{j=s+1}^p \delta_j^\sharp \leq (1+c_0)s^{1-1/q}\lVert \delta\rVert_q$ and hence \begin{align*}
\sum\limits_{j=s+1}^p \delta_j^\sharp \sqrt{\log(2p/j)}\leq \sqrt{\log(2p/s)}\sum\limits_{j=s+1}^p \delta_j^\sharp &\leq
(1+c_0)\lVert \delta\rVert_q s^{1-1/q}\sqrt{\log(2p/s)} \\
&\leq (1+c_0)\lVert \delta\rVert_q \left[ s \left(\sqrt{\log(2p/s)}\right)^{\frac{q}{q-1}}\right]^{1-1/q} \\
&\leq (1+c_0)\lVert \delta\rVert_q \left[ \sum\limits_{j=1}^s \left(\sqrt{\log(2p/j)}\right)^{\frac{q}{q-1}}\right]^{1-1/q}.
\end{align*}
In addition, the H\"{o}lder's inequality implies that \[
\sum\limits_{j=1}^s \delta_j^\sharp \sqrt{\log(2p/j)} \leq \lVert \delta\rVert_q \left[ \sum\limits_{j=1}^s \left(\sqrt{\log(2p/j)}\right)^{\frac{q}{q-1}}\right]^{1-1/q}.
\]
Combining these two statements, we obtain that $\delta\in \mathcal{C}_{WRE(q,s,1+c_0)}$, which leads to $\mathcal{C}_{SRE(q,s,c_0)}\subseteq\mathcal{C}_{WRE(q,s,1+c_0)}$. Therefore, the $WRE(q,s,1+c_0)$ condition implies the $SRE(q,s,c_0)$ condition. But the converse is not true, which means that there is no equivalence between these two conditions. A simple sufficient condition for $WRE(q,s,c_0)$ is given in Proposition 3 in Section 4. 

Now we are ready to state the main result for the Slope estimator, which acts as an extension of Theorem 6.1 in \cite{bellec2018}.

\begin{theorem}
	Let $s\in\{1,\cdots,p\}$, $\gamma\in(0,1)$ and $\tau\in[0, 1-\gamma)$. For any $2\leq q\leq \infty$, assume that the $WRE(q,s,c_0)$ condition holds with $c_0=c_0(\gamma,\tau)=\frac{1+\gamma+\tau}{1-\gamma-\tau}$. Let the tuning parameters $\{\lambda_j\}$ be defined by (\ref{weights}) with constant $A>(4+\sqrt{2})/\gamma$. Let $\delta_0\in(0,1)$. Then, on the event (\ref{event}), the Slope estimator $\hat{\beta}$ with tuning parameters $\{\lambda_j\}$ satisfies \begin{align}
	2\tau\lVert\hat{\beta}-\beta\rVert_{*}+\frac{1}{n}\lVert X\hat{\beta}-f\rVert_2^2\leq \frac{1}{n}\lVert X\beta-f\rVert_2^2+C'_{\gamma,\tau}(q,s,\delta_0)\Lambda_q^2(s)+4\sigma_s(\beta)_{*}, \label{slope1}
	\end{align}
	for all $\beta\in\mathbb{R}^p$ and all $f\in\mathbb{R}^n$, where \begin{align*}
	C'_{\gamma,\tau}(q,s,\delta_0):=(1+\gamma+\tau)^2\left(\frac{\log(1/\delta_0)}{s\log(2p/s)}\vee \frac{1}{\nu_q^2(s,c_0(\gamma,\tau))}\right).
	\end{align*}
	Furthermore, if $f=X\beta^{*}$ for some $\beta^{*}$, then on the event (\ref{event}), we have \begin{align}
	\lVert \hat{\beta}-\beta^{*}\rVert_{*}&\leq \frac{C'_{\gamma,\tau}(q,s,\delta_0)}{2\tau}\Lambda_q^2(s)+\frac{2}{\tau}\sigma_s(\beta^{*})_{*},  \label{slope2} \\ 
	\lVert \hat{\beta}-\beta^{*}\rVert_q&\leq \frac{C'_{\gamma,0}(q,s,\delta_0)}{1+\gamma}\Lambda_q(s),\quad (\text{when $\beta^{*}$ is $s$-sparse})  \label{slope3}  \\
    \lVert \hat{\beta}-\beta^{*}\rVert_q&\leq \left(\frac{2}{1+\gamma}\vee \frac{1-\gamma-\tau}{4\tau}\right)C'_{\gamma,\tau}(q,s,\delta_0)\Lambda_q(s)+\left(\frac{1}{1+\gamma}\vee \frac{1-\gamma-\tau}{\tau}\right)\Lambda_q(s)^{-1}\sigma_s(\beta^{*})_*. \label{slope4}
	\end{align}
 Here, note that since the tuning parameters $\{\lambda_j\}$ satisfy (\ref{weights}), so we have \begin{align}
	\Lambda_q(s)=\left(\sum_{j=1}^s \lambda_j^\frac{q}{q-1}\right)^{1-1/q}\leq s^{1/2-1/q}\left(\sum\limits_{j=1}^s \lambda_j ^2\right)^{1/2}\leq s^{1/2-1/q}\sqrt{\frac{A^2\sigma^2s\log(2ep/s)}{n}}=A\sigma s^{1-1/q}\sqrt{\frac{\log(2ep/s)}{n}}.
	\end{align}
\end{theorem}
\medskip

\noindent
{\bf Proof.} The proof follows the same route as the Proof of Theorem 1 with some careful modifications. Adopting (\ref{KKT}) with $h=2\lVert \cdot\rVert_{*}$, we have for all $\beta\in\mathbb{R}^p$ and all $f\in\mathbb{R}^n$, \begin{align}
2\tau\lVert \hat{\beta}-\beta\rVert_{*}+\frac{1}{n}\lVert X\hat{\beta}-f\rVert_2^2\leq \frac{1}{n}\lVert X\beta-f\rVert_2^2-\frac{1}{n}\lVert X(\hat{\beta}-\beta)\rVert_2^2+\Delta^{*}, \label{slopeKKT}
\end{align}
where $\Delta^{*}=2\tau \lVert \hat{\beta}-\beta\rVert_{*}+\frac{2}{n}\xi^{T}X(\hat{\beta}-\beta)+2\lVert \beta\rVert_{*}-2\lVert \hat{\beta}\rVert_{*}$. Let $u=\hat{\beta}-\beta$, then on the event (\ref{event}), Lemma 1 implies that \begin{align*}
\Delta^{*}&\leq 2(\tau\lVert \hat{\beta}-\beta\rVert_{*}+\lVert \beta\rVert_{*}-\lVert\hat{\beta}\rVert_{*})+2\max (\tilde{H}(u), G(u)) \\
&\leq 2\left((1+\tau)\lVert u\rVert_q \Lambda_q(s)-(1-\tau)\sum\limits_{j=s+1}^p \lambda_j u_j^{\sharp}+2\sigma_s(\beta)_{*}\right)+2\max (\tilde{H}(u), G(u)). 
\end{align*}

(i) If $\tilde{H}(u)\leq G(u)$, then the definitions of $\tilde{H}(u)$ and $G(u)$ yield \begin{align*}
\lVert u\rVert_q\leq \frac{G(u)}{\frac{4+\sqrt{2}}{\sqrt{n}}\sigma \left[\sum\limits_{j=1}^s\left(\sqrt{\log (2p/j)}\right)^{\frac{q}{q-1}}\right]^{1-1/q}}\leq \frac{G(u)}{\frac{4+\sqrt{2}}{\sqrt{n}}\sigma \sqrt{\sum\limits_{j=1}^s\log (2p/j)}}&\leq \frac{G(u)}{\frac{4+\sqrt{2}}{\sqrt{n}}\sigma \sqrt{s\log(2p/s)}}\\
&\leq \sqrt{\frac{\log(1/\delta_0)}{s\log(2p/s)}}\frac{\lVert Xu\rVert_2}{\sqrt{n}},
\end{align*}
where we use the facts that $q\in[2,\infty]$ and $\sum\limits_{j=1}^s \log(2p/j)\geq s\log(2p/s)$.
In addition, the setting of weights leads to \begin{align}
G(u)\leq \gamma \Lambda_2(s)\sqrt{\frac{\log(1/\delta_0)}{s\log(2p/s)}}\frac{\lVert Xu\rVert_2}{\sqrt{n}}\leq \gamma \Lambda_q(s)\sqrt{\frac{\log(1/\delta_0)}{s\log(2p/s)}}\frac{\lVert Xu\rVert_2}{\sqrt{n}},
\end{align}
since $\Lambda_2(s)\leq \Lambda_q(s)$ for any $q\in[2,\infty]$. Therefore, we obtain that \begin{align*}
\Delta^{*}&\leq 2(1+\tau)\lVert u\rVert_q\Lambda_q(s)+4\sigma_s(\beta)_{*}+2G(u) \\
&\leq 2(1+\gamma+\tau)\Lambda_q(s)\sqrt{\frac{\log(1/\delta_0)}{s\log(2p/s)}}\frac{\lVert Xu\rVert_2}{\sqrt{n}}+4\sigma_s(\beta)_{*} \\
&\leq (1+\gamma+\tau)^2\Lambda_q^2(s)\frac{\log(1/\delta_0)}{s\log(2p/s)}+\frac{\lVert Xu\rVert_2^2}{n}+4\sigma_s(\beta)_{*}.
\end{align*}

(ii) If $\tilde{H}(u)>G(u)$, in which case $\tilde{H}(u)\leq \gamma (\lVert u\rVert_q\Lambda_q(s)+\sum\limits_{j=s+1}^p \lambda_j u_j^{\sharp})$, then we have \begin{align*}
\Delta^{*}&\leq 2\left((1+\tau)\lVert u\rVert_q\Lambda_q(s)-(1-\tau)\sum\limits_{j=s+1}^p \lambda_j u_j^{\sharp}+2\sigma_s(\beta)_{*}\right)+2\tilde{H}(u) \\
&\leq \underbrace{2(1+\gamma+\tau)\lVert u\rVert_q\Lambda_q(s)-2(1-\gamma-\tau)\sum\limits_{j=s+1}^p\lambda_ju_j^{\sharp}}_{\Delta}+4\sigma_s(\beta)_{*}.
\end{align*}
If $\Delta\leq 0$, then the statement (\ref{slope1}) holds trivially in view of (\ref{slopeKKT}). If $\Delta>0$, then $\sum\limits_{j=s+1}^p \lambda_j u_j^{\sharp}\leq \frac{1+\gamma+\tau}{1-\gamma-\tau}\lVert u\rVert_q\Lambda_q(s)$. Furthermore, as a result of $\sum\limits_{j=1}^s \lambda_j u_j^{\sharp}\leq \lVert u\rVert_q\Lambda_q(s)$, it holds that \begin{align*}
\lVert u\rVert_{*}=\sum\limits_{j=1}^p \lambda_j u_j^\sharp\leq (1+\frac{1+\gamma+\tau}{1-\gamma-\tau})\lVert u\rVert_q\Lambda_q(s)=(1+c_0)\lVert u\rVert_q\Lambda_q(s).
\end{align*}
Hence, $u$ belongs to the cone $C_{WRE}(q,s,c_0)$. Then by using the $WRE(q,s,c_0)$ condition we have \begin{align*}
\Delta^{*}\leq \Delta+4\sigma_s(\beta)_{*}&\leq 2(1+\gamma+\tau)\lVert u\rVert_q\Lambda_q(s)+4\sigma_s(\beta)_{*} \\
&\leq \frac{2(1+\tau+\gamma)\Lambda_q(s)\lVert Xu\rVert_2/\sqrt{n}}{\nu_q(s,c_0)}+4\sigma_s(\beta)_{*}\\
&\leq \frac{(1+\tau+\gamma)^2\Lambda_q^2(s)}{\nu_q^2(s,c_0)}+\frac{\lVert Xu\rVert_2^2}{n}+4\sigma_s(\beta)_{*}.
\end{align*}

Consequently, for both cases it holds that $\Delta^{*}\leq \underbrace{(1+\gamma+\tau)^2\left(\frac{\log(1/\delta_0)}{s\log(2p/s)}\vee \frac{1}{\nu_q^2(s,c_0(\gamma,\tau))}\right)}_{C'_{\gamma,\tau}(q,s,c_0)} \Lambda_q^2(s)+\frac{\lVert Xu\rVert_2^2}{n}+4\sigma_s(\beta)_{*}$, which implies that \begin{align}
2\tau\lVert \hat{\beta}-\beta\rVert_{*}+\frac{1}{n}\lVert X\hat{\beta}-f\rVert_2^2\leq \frac{1}{n}\lVert X\beta-f\rVert_2^2+C'_{\gamma,\tau}(q,s,c_0)\Lambda_q^2(s)+4\sigma_s(\beta)_{*}
\end{align}
for all $\beta\in\mathbb{R}^p$ and all $f\in\mathbb{R}^n$. The proof of (\ref{slope1}) is completed. 

To prove (\ref{slope2}), let $f=X\beta^{*}$ and $\beta=\beta^{*}$ in (\ref{slope1}), we immediately obtain that \begin{align*}
2\tau\lVert \hat{\beta}-\beta^{*}\rVert_{*}\leq 2\tau\lVert \hat{\beta}-\beta^{*}\rVert_{*}+\frac{1}{n}\lVert X(\hat{\beta}-\beta^{*})\rVert_2^2\leq C'_{\gamma,\tau}(q,s,c_0)\Lambda_q^2(s)+4\sigma_s(\beta)_{*}.
\end{align*}

As for (\ref{slope3}), when $\beta^*$ is $s$-sparse (i.e., $\sigma_s(\beta^*)_*=0$), we set in what follows $u=\hat{\beta}-\beta^{*}$, $\tau=0$ so that $c_0=\frac{1+\gamma}{1-\gamma}$ and consider the same two Cases as above. 

(i) If $\tilde{H}(u)\leq G(u)$, then it is derived from \begin{align}
\frac{1}{n}\lVert X\hat{\beta}-f\rVert_2^2\leq \frac{1}{n}\lVert X\beta^{*}-f\rVert_2^2 -\frac{1}{n}\lVert X\hat{\beta}-f \rVert_2^2+\Delta^{*}, \label{errorbound}
\end{align} that \begin{align*}
\frac{1}{n}\lVert Xu\rVert_2^2=\frac{1}{n}\lVert X(\hat{\beta}-\beta^{*})\rVert_2^2\leq (1+\gamma)^2\Lambda_q^2(s)\left(\frac{\log(1/\delta_0)}{s\log(2p/s)}\right).
\end{align*}
Combining this inequality and the fact that \begin{align*}
\lVert u\rVert_q\leq \sqrt{\frac{\log(1/\delta_0)}{s\log(2p/s)}}\frac{\lVert Xu\rVert_2}{\sqrt{n}},
\end{align*}
we get \begin{align*}
\lVert u\rVert_q&\leq \sqrt{\frac{\log(1/\delta_0)}{s\log(2p/s)}}\cdot\sqrt{(1+\gamma)^2\Lambda_q^2(s)\left(\frac{\log(1/\delta_0)}{s\log(2p/s)}\right)} \\
&\leq \sqrt{\frac{C'_{\gamma,0}(q,s,\delta_0)}{(1+\gamma)^2}}\sqrt{(1+\gamma)^2\Lambda_q^2(s)\frac{C'_{\gamma,0}(q,s,\delta_0)}{(1+\gamma)^2}} =\frac{C'_{\gamma,0}(q,s,\delta_0)}{1+\gamma}\Lambda_q(s).
\end{align*}

(ii) If $\tilde{H}(u)>G(u)$, then $\Delta\geq \Delta^{*}$. It follows from (\ref{errorbound}) with $\beta=\beta^{*}$ and $f=X\beta^{*}$ that $\Delta^{*}\geq \frac{2}{n}\lVert Xu\rVert_2^2$. Thus $\Delta\geq 0$, which implies that $u$ belongs to the cone $C_{WRE}(q,s,\frac{1+\gamma}{1-\gamma})$. By applying the $WRE(q,s,\frac{1+\gamma}{1-\gamma})$ condition, we obtain that \begin{align*}
\lVert u\rVert_q\leq \frac{\lVert Xu\rVert_2}{\sqrt{n}\nu_q(s,\frac{1+\gamma}{1-\gamma})}.
\end{align*}
Combining this bound and the fact that \begin{align*}
\frac{1}{n}\lVert Xu\rVert_2^2\leq \frac{(1+\gamma)^2\Lambda_q^2(s)}{\nu_q^2(s,\frac{1+\gamma}{1-\gamma})},
\end{align*}
we have \begin{align*}
\lVert u\rVert_q\leq \frac{\lVert Xu\rVert_2}{\sqrt{n}\nu_q(s,\frac{1+\gamma}{1-\gamma})}&\leq \frac{1}{\nu_q(s,\frac{1+\gamma}{1-\gamma})}\sqrt{\frac{(1+\gamma)^2\Lambda_q^2(s)}{\nu_q^2(s,\frac{1+\gamma}{1-\gamma})}}\leq \frac{\sqrt{C'_{\gamma,0}(q,s,\delta_0)}}{1+\gamma}\sqrt{C'_{\gamma,0}(q,s,\delta_0)\Lambda_q^2(s)}=\frac{C'_{\gamma,0}(q,s,\delta_0)}{1+\gamma}\Lambda_q(s),
\end{align*}
which completes the proof of (\ref{slope3}).

Finally, regarding the result (\ref{slope4}) that allows $\beta^{*}$ to be not exactly sparse,  we let $f=X\beta^{*}$, $\beta=\beta^{*}$ and $u=\hat{\beta}-\beta^{*}$ in (\ref{slope1}) and obtain \begin{align}
 2\tau\lVert u\rVert_*+\frac{1}{n}\lVert Xu\rVert_2^2\leq C'_{\gamma,\tau}(q,s,\delta_0)\Lambda_q^2(s)+4\sigma_s(\beta^{*})_*.
\end{align}
Then, if $u\in\mathcal{C}_{WRE}(q,s,c_0=c_0(\gamma,\tau)=\frac{1+\gamma+\tau}{1-\gamma-\tau})$, the $WRE(q,s,c_0)$ condition implies that \begin{align}
\lVert u \rVert_q\leq \frac{\lVert Xu\rVert_2/\sqrt{n}}{\nu_q(s,c_0(\gamma,\tau))}&\leq \frac{1}{\nu_q(s,c_0(\gamma,\tau))}\sqrt{C'_{\gamma,\tau}(q,s,\delta_0)\Lambda_q^2(s)+4\sigma_s(\beta)_*} \nonumber \\
& \leq \frac{\sqrt{C'_{\gamma,\tau}(q,s,\delta_0)}}{1+\gamma}\sqrt{C'_{\gamma,\tau}(q,s,\delta_0)\Lambda_q^2(s)+4\sigma_s(\beta)_*} \nonumber \\
&\leq\frac{C'_{\gamma,\tau}(q,s,\delta_0)}{1+\gamma}\Lambda_q(s)+\frac{1}{1+\gamma}\left(2\sqrt{C'_{\gamma,\tau}(q,s,\delta_0)\Lambda_q(s)}\cdot\Lambda_q(s)^{-1/2}\sqrt{\sigma_s(\beta^*)_*}\right) \nonumber \\
&\leq \frac{2}{1+\gamma}C'_{\gamma,\tau}(q,s,\delta_0)\Lambda_q(s)+\frac{1}{1+\gamma}\Lambda_q(s)^{-1}\sigma_s(\beta^*)_*.
\label{slope_compressible1}
\end{align}
Otherwise, if $u\notin\mathcal{C}_{WRE}(q,s,c_0=c_0(\gamma,\tau)=\frac{1+\gamma+\tau}{1-\gamma-\tau})$, then by using (\ref{slope2}) we have \begin{align}
\lVert u\rVert_q<\frac{1}{1+c_0}\Lambda_q(s)^{-1}\lVert u\rVert_*&<\frac{1-\gamma-\tau}{2}\Lambda_q(s)^{-1}\left(\frac{C'_{\gamma,\tau}(q,s,\delta_0)}{2\tau}\Lambda_q^2(s)+\frac{2}{\tau}\sigma_s(\beta^{*})_{*}\right) \nonumber \\
&<\frac{1-\gamma-\tau}{4\tau}C'_{\gamma,\tau}(q,s,\delta_0)\Lambda_q(s)+\frac{1-\gamma-\tau}{\tau}\Lambda_q(s)^{-1}\sigma_s(\beta^*)_*.\label{slope_compressible2}
\end{align}
As a consequence, (\ref{slope4}) holds due to (\ref{slope_compressible1}) and (\ref{slope_compressible2}). \\
	
With $\gamma=1/2$, $\tau=1/4$ and $\delta_0=(\frac{s}{2p})^{\frac{s}{\nu_q^2(s,7)}}$, it is straightforward to get the following corollary.

\begin{cor}
		Let $s\in\{1,\cdots,p\}$. For any $2\leq q\leq \infty$, assume that the $WRE(q,s,7)$ condition holds. Let the tuning parameters $\{\lambda_j\}$ be defined by (\ref{weights}) with constant $A>2(4+\sqrt{2})$. Then, with probability at least $1-\frac{1}{2}(\frac{s}{2p})^{\frac{s}{\nu_q^2(s,7)}}$, the Slope estimator $\hat{\beta}$ with tuning parameters $\{\lambda_j\}$ satisfies \begin{align*}
	\frac{1}{2}\lVert\hat{\beta}-\beta\rVert_{*}+\frac{1}{n}\lVert X\hat{\beta}-f\rVert_2^2\leq \frac{1}{n}\lVert X\beta-f\rVert_2^2+\frac{49\Lambda_q^2(s)}{16\nu_q^2(s,7)}+4\sigma_s(\beta)_{*}, 
	\end{align*}
	for all $\beta\in\mathbb{R}^p$ and all $f\in\mathbb{R}^n$. Furthermore, if $f=X\beta^{*}$ for some $\beta^{*}$, then with probability at least $1-\frac{1}{2}(\frac{s}{2p})^{\frac{s}{\nu_q^2(s,7)}}$, we have \begin{align*}
	\lVert \hat{\beta}-\beta^{*}\rVert_{*}&\leq \frac{49\Lambda_q^2(s)}{8\nu_q^2(s,7)}+8\sigma_s(\beta^{*})_{*},   \\ 
	\lVert \hat{\beta}-\beta^{*}\rVert_q&\leq\frac{49\Lambda_q(s)}{24\nu_q^2(s,7)},\quad (\text{when $\beta^{*}$ is $s$-sparse})    \\
	\lVert \hat{\beta}-\beta^{*}\rVert_q&\leq \frac{49\Lambda_q(s)}{12\nu_q^2(s,7)}+\Lambda_q(s)^{-1}\sigma_s(\beta^{*})_*. 
	\end{align*}
\end{cor}

\section{Matching Lower Bounds}

In this section, we move on to discuss the lower bounds for the $\ell_q$ estimation error over both the $\ell_0$-ball $B_0(s):=\{\beta\in\mathbb{R}^p:\lVert \beta\rVert_{0}\leq s\}$ and the $\ell_r$-balls $B_r(R_r):=\{\beta\in\mathbb{R}^p:\lVert \beta\rVert_{r}^r=\sum_{j=1}^p |\beta_j|^r\leq R_r\}$ with $0<r<1$.

\subsection{Over the $\ell_0$-ball}
In what follows, we borrow the main result from Section 7 in \cite{bellec2018}. Let us introduce the maximal $s$-sparse eigenvalue of the design matrix $X$ defined as \begin{align}
\bar{\theta}_{\mathrm{max}}(X,s)=\max\limits_{\delta\in B_0(s)\setminus\{0\}} \frac{\lVert X\delta\rVert_2}{\sqrt{n}\lVert \delta\rVert_2},
\end{align}
with $\bar{\theta}_{\mathrm{max}}(X,1)=\max_{j=1,\cdots,p}\lVert Xe_j\rVert_2/\sqrt{n}$, and $\psi_{n,q}=\sigma s^{1/q}\sqrt{\frac{\log(ep/s)}{n}}$ for any $1\leq q\leq \infty$. Then, according to Theorem 7.1 of \cite{bellec2018}, we have the following lemma, which establishes the minimax lower bound for the $\ell_q$ estimation error over the $\ell_0$-ball $B_0(s)$ for any design matrix $X$.

\begin{lemma}
	Let $p\geq 2$, $s\in[1,p/2]$, $n\geq 1$ be integers and $1\leq q\leq \infty$. Assume that $y=X\beta^{*}+\xi$ with $\xi\in\mathbb{R}^n\sim N(0,\sigma^2 I_{n\times n})$, $\sigma>0$. Then it holds that there exist positive constants $\bar{b},\bar{c}$ depending on $q$ such that \begin{align}
	\inf\limits_{\hat{\tau}}\inf\limits_{X}\sup\limits_{\beta^{*}\in B_0(s)} P\Big(\lVert \hat{\tau}-\beta^{*}\rVert_q\geq \psi_{n,q}\bar{b}^{-1}\bar{\theta}_{\mathrm{max}}^{-1}(X,1)\Big)\geq \bar{c},
	\end{align}
	where $\inf_{\hat{\tau}}$ denotes the infimum over all estimators $\hat{\tau}$ of $\beta^{*}$ and $\inf_{X}$ denotes the infimum over all matrices $X\in\mathbb{R}^{n\times p}$.
\end{lemma}

Hence, when the true parameter vector $\beta^{*}$ is exactly $s$-sparse, the upper bounds for the Lasso and Slope estimators obtained in the previous sections show that both the Lasso estimator with a tuning parameter of order $\sqrt{\log(p/s)/n}$ and the Slope estimator with weights defined by (\ref{weights}) achieve the minimax optimal rate $\sqrt{s\log(p/s)/n}$ for the $\ell_2$ estimation error. In fact, the Lasso estimator with a tuning parameter of order $\sqrt{\log(p/s)/n}$ achieves the minimax optimal rate $s^{1/q}\sqrt{\log(p/s)/n}$ for the $\ell_q$ estimation error with $1\leq q\leq 2$, while a near-optimal rate $s^{1-1/q}\sqrt{\log(p/s)/n}$ for the $\ell_q$ estimation error with $2<q\leq \infty$. The Slope estimator with weights defined by (\ref{weights}) achieves the near-optimal rate $s^{1-1/q}\sqrt{\log(p/s)/n}$ for the $\ell_q$ estimation error with $2<q\leq \infty$ as well.

\subsection{Over the $\ell_r$-balls}
As shown in the following proposition, when the true parameter vector $\beta^{*}$ is not exactly sparse, the error of its best $s$-term approximation with respect to $\lVert \cdot\rVert_{*}$ can be well controlled by its $\ell_r$ norms with $0<r<1$ so that the non-convex $\ell_r$-balls $B_r(R_r)$ serve as good models for vectors that are not exactly sparse. 

\begin{proposition}
	For any $s\in \{1,\cdots,p\}$, $\lambda_1\geq \lambda_2\geq \cdots\geq\lambda_p\geq 0$, $0<r<1$ and $\beta\in\mathbb{R}^p$, \begin{align}
	\sigma_s(\beta)_{*}\leq \lambda_s \cdot s^{\frac{r}{1-r}}\lVert \beta\rVert_r.
	\end{align}
	In particular, if $\lambda_1=\lambda_2=\cdots=\lambda_p=\lambda$ for some $\lambda>0$, then $\sigma_s(\beta)_{*}=\lambda\sigma_s(\beta)_1$, and we have \begin{align}
	\sigma_s(\beta)_1\leq s^{\frac{r}{1-r}}\lVert \beta\rVert_r.
	\end{align}
\end{proposition}

\noindent
{\bf Proof.} Based on the definition of $\sigma_s(\beta)_{*}$,  we have \begin{align*}
\sigma_s(\beta)_{*}=\sum\limits_{j=s+1}^p \lambda_j \beta_j^{\sharp}&\leq (\beta_s^{\sharp})^{1-r}\sum\limits_{j=s+1}^p \lambda_j(\beta_j^{\sharp})^r \\
&\leq \left(\frac{1}{s}\sum\limits_{j=1}^s (\beta_j^{\sharp})^r\right)^{\frac{1-r}{r}}\cdot\sum\limits_{j=s+1}^p \lambda_j(\beta_j^{\sharp})^r \\
&\leq  \left(\frac{1}{s}\lVert \beta\rVert_r^r\right)^{\frac{1-r}{r}}\cdot\lambda_s\lVert \beta\rVert_r^r \\
&\leq  \lambda_s \cdot s^{\frac{r}{1-r}}\lVert \beta\rVert_r,
\end{align*}
which completes the proof. \\

Unfortunately, how to generalize Lemma 4 mentioned above and thus obtain the non-asymptotic minimax lower bound for the $\ell_q$ estimation error with general $1\leq q\leq \infty$ and for any design matrix $X$ over the $\ell_r$-balls remains open. The corresponding discussion is left for future work.

\section{Restricted Eigenvalue Conditions}

In this section, we study the relationships between the proposed RE conditions. 

\begin{proposition}
Let $c_0>0$ and $s\in\{1,\cdots,p\}$. For any $q\in[2,\infty]$ we have
if condition $SRE(s^{\frac{2(q-1)}{q}},c_0)$ holds, then condition $SRE(q,s,c_0)$ holds and $\theta_q(s,c_0)\geq \theta(s^{\frac{2(q-1)}{q}},c_0)$.
\end{proposition}

\noindent
{\bf Proof.} For any $q\in[2,\infty]$, the cone $\mathcal{C}_{SRE(q,s,c_0)}$ of the $SRE(q,s,c_0)$ is the set of all non-zero $\delta\in\mathbb{R}^p$ such that $\lVert \delta\rVert_1\leq (1+c_0) s^{1-1/q}\lVert \delta\rVert_q$. By using the fact that $\lVert \delta\rVert_q\leq \lVert \delta\rVert_2$, we have $\lVert \delta\rVert_1\leq (1+c_0) \sqrt{s^{2-2/q}}\lVert \delta\rVert_2$, so that $\mathcal{C}_{SRE(q,s,c_0)}\subseteq \mathcal{C}_{SRE(s^{2-2/q},c_0)}$. As a consequence,
\begin{align*}
\theta_q(s,c_0):=\min\limits_{\delta\neq 0, \delta\in\mathcal{C}_{SRE(q,s,c_0)}} \frac{\lVert X\delta\rVert_2}{\sqrt{n}\lVert \delta\rVert_q}\geq \min\limits_{\delta\neq 0, \delta\in\mathcal{C}_{SRE(s^{2-2/q},c_0)}} \frac{\lVert X\delta\rVert_2}{\sqrt{n}\lVert \delta\rVert_2}\cdot \frac{\lVert \delta\rVert_2}{\lVert \delta\rVert_q} 
\geq \theta(s^{\frac{2(q-1)}{q}},c_0). 
\end{align*} 

\medskip
Hence, based on this proposition and the results in \cite{bellec2018} for the $SRE(s,c_0)$ condition, we can conclude that for a large class of random matrices with i.i.d. rows, the $SRE(q,s,c_0)$ condition holds with high probability when $s^{\frac{2(q-1)}{q}}\log(ep/s)\leq cn$ with $c>0$ being a constant. In addition, by using the following proposition we are able to obtain that the $WRE(q,s,c_0)$ holds with high probability whenever $s^{\frac{2(q-1)}{q}}\log^2(ep/s)\leq Cn$ with $C>0$ being some constant.

\begin{proposition}
	Let $c_0>0$, $s\in\{1,\cdots,p\}$ and $s_q=\lceil s\left(\sqrt{\log(2ep/s)/\log 2}\right)^{\frac{q}{q-1}} \rceil $. For any $q\in[2,\infty]$ we have if condition $SRE(q,s_q,c_0)$ holds, then condition $WRE(q,s,c_0)$ holds and $\nu_q(s,c_0)\geq \theta_q(s_q,c_0)$.
\end{proposition}

\noindent
{\bf Proof.} If $\delta\in \mathcal{C}_{WRE(q,s,c_0)}$, then\begin{align*}
(1+c_0)\left(\sum\limits_{j=1}^s \lambda_j ^{\frac{q}{q-1}}\right)^{1-1/q}\lVert\delta\rVert_q\geq \lVert \delta\rVert_{*}=\sum\limits_{j=1}^p\lambda_j\delta_j^\sharp\geq \lambda_p \sum\limits_{j=1}^p \delta_j^\sharp=\lambda_p\lVert \delta\rVert_1=\frac{A\sigma}{\sqrt{n}} \sqrt{\log{2}}\lVert \delta\rVert_1.
\end{align*}
Meanwhile, $\left(\sum\limits_{j=1}^s \lambda_j ^{\frac{q}{q-1}}\right)^{1-1/q}\leq s^{1/2-1/q}\left(\sum\limits_{j=1}^s \lambda_j ^2\right)^{1/2}\leq\frac{A\sigma}{\sqrt{n}}  s^{1/2-1/q}\sqrt{s\log(2ep/s)}=\frac{A\sigma}{\sqrt{n}} \left[s\left(\sqrt{\log(2ep/s)}\right)^{\frac{q}{q-1}}\right]^{1-1/q}$. Thus, we have $\lVert \delta\rVert_1\leq (1+c_0)\left[s\left(\sqrt{\log(2ep/s)/\log 2}\right)^{\frac{q}{q-1}}\right]^{1-1/q}\lVert \delta\rVert_q$ and $\delta\in\mathcal{C}_{SRE(q,s_q,c_0)}$, which implies that $\mathcal{C}_{WRE(q,s,c_0)}\subseteq\mathcal{C}_{SRE(q,s_q,c_0)}$ and completes the proof.

\section{Conclusion}

In this paper, we have generalized the sharp results for the Lasso and the Slope estimators in \cite{bellec2018} to allow the case that the true parameter vector is not exactly sparse, and obtained the results for any $q\in[1,\infty]$ by using the corresponding extended Restricted Eigenvalue conditions. What we have done could be further developed to other Slope type estimators including the Square-Root Slope \cite{derumigny2018}, the Group Slope \cite{brzyski2019} and the Sorted Concave Penalized estimator \cite{feng2019}, which are left for future works.



\begin{thebibliography}{}
	
\bibitem{bellec2018}	
Bellec, P.C., Lecué, G., \& Tsybakov, A.B. (2018). Slope meets lasso: improved oracle bounds and optimality. \textit{ Annals of Statistics}, 46(6B), 3603-3642.
	
\bibitem{bickel2009}	
Bickel, P.J., Ritov, Y.A., \& Tsybakov, A.B. (2009). Simultaneous analysis of Lasso and Dantzig selector. \textit{Annals of Statistics}, 37(4), 1705-1732.

\bibitem{bogdan2015}
Bogdan, M., Van Den Berg, E., Sabatti, C., Su, W., \& Candès, E.J. (2015). SLOPE—adaptive variable selection via convex optimization. \textit{Annals of Applied Statistics}, 9(3), 1103.

\bibitem{brzyski2019}
Brzyski, D., Gossmann, A., Su, W., \& Bogdan, M. (2019). Group slope–adaptive selection of groups of predictors. \textit{Journal of the American Statistical Association}, 114(525), 419-433.

\bibitem{buhlmann2011}
Bühlmann, P., \& Van De Geer, S. (2011). Statistics for high-dimensional data: methods, theory and applications. Springer Science \& Business Media.

\bibitem{derumigny2018}
Derumigny, A. (2018). Improved bounds for square-root lasso and square-root slope. \textit{Electronic Journal of Statistics}, 12(1), 741-766.

\bibitem{feng2019}
Feng, L., \& Zhang, C. H. (2019). Sorted concave penalized regression. \textit{Annals of Statistics}, 47(6), 3069-3098.

\bibitem{lecue2018}
Lecué, G., \& Mendelson, S. (2018). Regularization and the small-ball method i: sparse recovery. \textit{Annals of Statistics}, 46(2), 611-641.

\bibitem{raskutti2011}
Raskutti, G., Wainwright, M.J., \& Yu, B. (2011). Minimax rates of estimation for high-dimensional linear regression over $\ell_q $-balls. \textit{IEEE Transactions on Information Theory}, 57(10), 6976-6994.

\bibitem{su2016}
Su, W., \& Candes, E. (2016). SLOPE is adaptive to unknown sparsity and asymptotically minimax. \textit{Annals of Statistics}, 44(3), 1038-1068.

\bibitem{tang2011}
Tang, G., \& Nehorai, A. (2011). Performance analysis of sparse recovery based on constrained minimal singular values. \textit{IEEE Transactions on Signal Processing}, 59(12), 5734-5745.

\bibitem{tibshirani1996}
Tibshirani, R. (1996). Regression shrinkage and selection via the lasso. \textit{Journal of the Royal Statistical Society: Series B (Methodological)}, 58(1), 267-288.

\bibitem{zhou2019a}
Zhou, Z., \& Yu, J. (2019). Sparse recovery based on $q$-ratio constrained minimal singular values. \textit{Signal Processing}, 155, 247-258.

\bibitem{zhou2019b}
Zhou, Z., \& Yu, J. (2019). On $q$-ratio CMSV for sparse recovery. \textit{Signal Processing}, 165, 128-132.
\end{thebibliography}
\end{document}